\documentclass[12pt]{article}
\newtheorem{theorem}{Theorem}

\newenvironment{proof}{{\bf Proof.}}{\hfill\rule{2mm}{2mm}}

\newtheorem{remarka}{Remark}

\def\lcs#1{{\rm lcs}(#1)}
\usepackage{array}
\newlength{\cellwid}
\makeatletter
\newenvironment{latinsq}[1][00]{%
 \def\centcol{\centering \let\\=\tabularnewline \CellStrut}
 \arraycolsep0pt
 \setbox\@tempboxa\hbox{#1}%
 \cellwid\ht\@tempboxa  \advance\cellwid\dp\strutbox
 \ifdim\cellwid<\wd\@tempboxa 
   \@tempdima .5\wd\@tempboxa \advance\@tempdima -.5\cellwid 
   \advance\cellwid \@tempdima \advance\@tempdima\dp\strutbox
 \else
   \@tempdima\dp\strutbox
 \fi
 \edef\CellStrut{\vrule
   width\z@ height\the\cellwid depth\the\@tempdima \relax}
 \advance\cellwid\@tempdima \advance\cellwid-2\arraycolsep
 \array{|>{\centcol}p{\cellwid}|*{20}{>{\centcol}p{\cellwid}|}}}%
 {\endarray}

\title{\bf A lower bound for the size of the largest critical sets in
Latin squares\footnote{This research in part supported by a grant
from IPM (No. 81050022).}}
\author{
{\bf Hamed Hatami$^a$ and
Ebadollah S. Mahmoodian$^b$} \\
\\[1mm]
Institute for Studies on Theoretical Physics and Mathematics (IPM) \\
{and}
\\
$^a${\small\it Department of Computer Engineering} \\
$^b${\small\it Department of Mathematical Sciences}\\
Sharif University of Technology \\
P.O. Box 11365--9415, Tehran, I.R. Iran}
\date{}
\begin{document}
\maketitle
\begin{abstract}
A {\sf critical set } in an $n \times n$ array is a set $C$ of
given entries, such that there exists a unique extension of  $C$
to  an $n\times n$ Latin square and no proper subset of $C$ has
this property. The cardinality of the largest critical set in any
Latin square of order $n$ is denoted by $\lcs{n}$.  We give a
lower bound for $\lcs{n}$ by showing that
 $\lcs{n} \geq n^2(1-\frac{2 + \ln 2}{\ln n})
 +n(1+\frac {\ln \left( 8 \pi \right)} {\ln n})-\frac{\ln 2}{\ln n}.$
\end{abstract}

\section{Introduction}        
A {\sf Latin square} of order $n$ is an $n$ $\times$ $n$ array of
integers chosen from the set $X = \{1,2, \ldots, n\}$ such that
each element of $X$ occurs exactly once in each row and exactly
once in each column. A Latin square can also be written as a set
of ordered triples $\{ (i,j;k) \mid$  symbol $k$ occurs in cell
$(i,j)$ of the array$\}$.

A {\sf partial Latin square} $P$ of order $n$ is an $n\times n$
array with entries chosen from the set $X = \{1,2, \ldots, n\}$,
such that each element of $X$ occurs at most  once in each row
and at most once in each column. Hence there are cells in the
array that may be empty, but the cells that are filled have been
filled so as to conform with the Latin property of the array. Let
$P$ be a partial Latin square of order $n$. Then $| P |$ is said
to be the {\sf size} of the partial Latin square and the set of
cells ${\cal S}_P=\{(i,j) \mid (i,j;k)\in P\}$ is said to
determine the {\sf shape} of $P$.

A partial Latin square $C$ contained in a Latin square $L$ is said
to be {\sf uniquely completable} if $L$ is the only Latin square
of order $n$ with $k$ in the cell $(i,j)$ for every $(i,j;k) \in
C$. A {\sf critical set} $C$ contained in a Latin square $L$ is a
partial Latin square that is uniquely completable and no proper
subset of $C$ satisfies this requirement. The name ``critical
set'' and the concept were invented by a statistician, John
Nelder, about 1977, and his ideas were first published in a
note~{\bf\cite{Nelder77}}.  This note posed the problem of giving
a formula for the size of the largest and smallest critical sets
for a Latin square of a given order. Let $ \lcs{n}$ denote the
size of the {\sf largest critical set} in any Latin square of
order $n$.  Nelder~{\bf\cite{Nelder?}} constructed a critical set
of size   $(n^2-n)/2$ for the $n \times n$ back circulant Latin
square. He conjectured that $ \lcs{n} =(n^2-n)/2$. This equality
was shown to be false in 1978, when Curran and van
Rees~{\bf\cite{MR80j:05022}}, found that $\lcs{4} \geq 7$. The
following is  an example of a largest critical set of size $11$
for a $5\times 5$ Latin square, taken
from~{\bf\cite{BeanMahmoodian}}, which also contradicts Nelder's
conjecture.

$$
\begin{latinsq}[1]
\hline  2&&4&3& \\
\hline&&1&2& \\
\hline&2&3&1& \\
\hline  3&1&2&& \\
\hline&&&& \\
\hline
\end{latinsq}
$$
\noindent

In the following table some known values of $\lcs{n}$ for $n \leq
6$ are listed,
$$
\begin{array}{c|cccccc}
n           & 1 & 2 & 3 & 4 &  5 &  6 \\ \hline
$\lcs{$n$}$ & 0 & 1 & 3 & 7 & 11 & 18  \\
\end{array}
$$
and in the following table some known lower bounds for $\lcs{n}$
 are shown for $7 \le n \le 10$,
$$
\begin{array}{c|cccc}
n             &  7 &  8 &  9 & 10  \\ \hline
\lcs{$n$} \ge & 25 & 37 & 44 & 57   \\
\end{array}
$$
See~{\bf\cite{BeanMahmoodian}} for the references. Recently Bean
and Mahmoodian~{\bf\cite{BeanMahmoodian}} have found the upper
bound $\lcs{n} \leq n^2-3n+3$. Nelder's $(n^2-n)/2$ is the best
lower bound that is found for $\lcs{n}$ so far. In this note we
improve this bound asymptotically for $n$ large enough ($n \ge
195$).


\section{A lower bound for $\lcs{n}$}

\begin{theorem}
For any integer $n$ we have,
\[
\lcs{n} \ge n^2(1-\frac{2 + \ln 2}{\ln n}) +n(1+\frac {\ln \left(
8 \pi \right)} {\ln n})-\frac{\ln 2}{\ln n}.
\]
\end{theorem}
\begin{proof} By Theorem 17.2 in~{\bf\cite{vanLintWilson}},
as a result of van der Warden conjecture, we know the following
bound for $L(n)$, the number of Latin squares of order~$n$: \ \
$L(n) \ge \frac{(n!)^{2n}}{n^{n^2}} $.

If in a partial Latin square all the entries, except the entries
of the first row and the first column be given, then it is
uniquely completable. So every Latin square has at least one
critical set which has no intersection with its first row and
first column. And also obviously the number of these critical sets
is greater than or equal to $L(n)$. For choosing the shape of such
a critical set we have at most $2^{(n-1)^2}$ ways, and for
choosing the entries of each given shape we have at most
$n^{\lcs{n}}$ different ways. So the number of critical sets is
less than or equal to $2^{n^2-2n+1}n^{\lcs{n}}$. Thus the
following inequalities hold:
\[ \frac{(n!)^{2n}}{n^{n^2}} \leq L(n) \leq 2^{n^2-2n+1}n^{\lcs{n}}. \]

Now by Stirling's approximation formula, (see for example
{\bf\cite{CLR}}), we can replace $n!$ with a smaller value
$\sqrt{2\pi n} \left( n \over e \right) ^ n$. So
\[ \frac{(2
\pi)^n n^{2n^2+n}}{e^{2n^2} n^{n^2}} \leq
2^{n^2-2n+1}n^{\lcs{n}},\] or
\[ \frac{(2 \pi)^n n^{n^2+n}}{e^{2n^2} 2^{n^2-2n+1}}
\leq n^{\lcs{n}}. \] Thus, \ \ $ n \ln (2 \pi) + (n^2+n) \ln n -
2n^2 - (n^2-2n+1) \ln 2 \leq \lcs{n}\ln n. $ \ This implies that $
n^2(1-\frac{2+\ln 2}{\ln n})+n(1+\frac {2{\ln 2}+\ln \left( 2 \pi
\right)} {\ln n})-\frac{\ln 2}{\ln n} \leq \lcs{n}.$
\end{proof}

\noindent {\bf Note.} Stinson and van
Rees~{\bf\cite{MR84g:05036}} have shown that
 $\lcs{2^m} \geq 4^m-3^m$. This lower bound for $n=2^m$, is better
 than the bound given in Theorem~1.
 \newpage
\def\cprime{$'$}

\end{document}